\documentclass[11pt]{amsart}

\usepackage{amssymb,amsfonts,amsmath,amsthm}
\usepackage[all,arc]{xy}
\usepackage{enumerate}
\usepackage{mathrsfs}
\usepackage[toc,page]{appendix}
\usepackage[left=3cm, right=3cm, bottom=3cm]{geometry}
\usepackage{tabularx}
\usepackage{graphicx}
\usepackage{url}
\usepackage{color}
\usepackage{tikz-cd}  
\usepackage{mathdots} 
\usepackage{romannum} 
\usepackage{hyperref} 

\newtheorem{thm}{Theorem}[section]

\newtheorem*{thm*}{Theorem}
\newtheorem*{cor*}{Corollary}
\newtheorem*{prop*}{Proposition}
\newtheorem{cor}[thm]{Corollary}
\newtheorem{prop}[thm]{Proposition}
\newtheorem{lem}[thm]{Lemma}
\newtheorem{conj}[thm]{Conjecture}

\theoremstyle{definition}

\newtheorem{exmp}[thm]{Example}

\newtheorem*{notn*}{Notation}

\theoremstyle{remark}
\newtheorem{rem}[thm]{Remark}

\newtheorem*{idea*}{Idea}

\newcommand{\Spec}{{\rm Spec}}

\newcommand{\sshf}[1]{\mathcal{O}_{#1}} 
\newcommand{\shf}[1]{\mathscr{#1}} 
\newcommand{\prj}[1]{\mathbb{P}^{#1}} 

\newcommand{\iso}{\simeq} 

\makeatletter
\let\c@equation\c@thm
\makeatother
\numberwithin{thm}{section}
\numberwithin{equation}{section}

\title[On the Image of the Hitchin Morphism over Algebraic Surfaces]{On the Image of Hitchin Morphism for Algebraic Surfaces: the Case ${\rm GL}_n$}

\author{Lei Song and Hao Sun}

\begin{document}
\pagenumbering{arabic}
\maketitle
\begin{abstract}
The Hitchin morphism is a map from the moduli space of Higgs bundles $\mathscr{M}_X$ to the Hitchin base $\mathscr{B}_X$, where $X$ is a smooth projective variety. When $X$ has dimension at least two, this morphism is not surjective in general. Recently, Chen-Ng\^o introduced a closed subscheme $\mathscr{A}_X$ of $\mathscr{B}_X$, which is called the space of spectral data. They proved that the Hitchin morphism factors through $\mathscr{A}_X$ and conjectured that $\mathscr{A}_X$ is the image of the Hitchin morphism. We prove that when $X$ is a smooth projective surface, this conjecture is true for vector bundles. Moreover, we show that $\mathscr{A}_X$, for any dimension, is invariant under any proper birational morphism, and apply the result to study $\mathscr{A}_X$ for ruled surfaces. 
\end{abstract}

\flushbottom



\renewcommand{\thefootnote}{\fnsymbol{footnote}}
\footnotetext[1]{MSC2010 Class: 14D20, 14J60}

\section{Introduction}
Throughout the paper, we work over an algebraically closed field $k$ of characteristic zero. Let $\mathcal{M}_X$ be the moduli space of semistable Higgs bundles of rank $n$ over a smooth projective variety $X$ over $k$. Let
\begin{align*}
    \mathscr{B}_X= \bigoplus\limits_{i=1}^n H^0(X, S^n T^*X)
\end{align*}
be the Hitchin base, where $T^*X$ is the cotangent bundle of $X$. The Hitchin morphism $h_X: \mathcal{M}_X \rightarrow \mathscr{B}_X$ was introduced in Hitchin's seminal work \cite{Hit1987} for algebraic curves. The morphism was proved to be dominant by Beauville-Narasimhan-Ramanan \cite{BNR}, and was later on proved to be proper by Nisture \cite{Nit} and Simpson \cite{Simp2}, and hence surjective in the case of algebraic curves. The properness and surjectivity of the Hitchin map are elementary and important properties in the study of $\mathcal{M}_X$, which includes for instance connected components \cite{Goth,KSZ211}, cohomology ring \cite{HauTha1,HauTha2}, Langlands duality \cite{ChenZhu2,DoPan} and mirror symmetry \cite{HauThaMLH}.

In a more general setting and for higher dimensional varieties, Simpson showed that the Hitchin morphism is still proper \cite{Simp1, Simp2}, but the surjectivity does not hold in general. Recently, Chen and Ng\^o took an attempt to understand the image of the Hitchin morphism and obtained a higher dimensional analogue of the BNR correspondence \cite{BNR}. In \cite{CN,CN2}, they considered the moduli stack $\shf{M}_X$ of Higgs bundles on $X$, and introduced a closed subscheme $\mathfrak{i}_X: \mathscr{A}_X \hookrightarrow \mathscr{B}_X$, where
\begin{align*}
    \mathscr{A}_X={\rm Sect}(X,{\rm Chow}^n(T^*_X/X))
\end{align*}
is called \emph{the space of spectral data} (see \S 2). They showed the Hitchin morphism factors through the space of spectral data
\begin{center}
\begin{tikzcd}
            & \shf{M}_X \arrow[d,"h_X"] \arrow[ld, dotted, "sd_X" description]  \\
            \shf{A}_X \arrow[r, hook, "\mathfrak{i}_X"] & \shf{B}_X,
\end{tikzcd}
\end{center}
and conjectured that the morphism $sd_X$ is surjective. The morphism $sd_X$ is called \emph{the spectral data morphism}.

\begin{conj}[Conjecture 5.2 in \cite{CN2}]
For every point $a \in \shf{A}_X(k)$, the fiber $sd^{-1}_X(a)$ is nonempty.
\end{conj}
In fact, the conjecture is stated for $G$-Higgs bundles, where $G$ is a split reductive group. In case the dimension $d=2$ and $G={\rm GL}_n$, the conjecture is verified by Chen and Ng\^o for spectral data in an open subset $\shf{A}^{\heartsuit}_X(k) \subseteq \mathscr{A}_X(k)$; and it is verified for some minimal surfaces, such as ruled surfaces and nonisotrivial elliptic surfaces for all $a\in \shf{A}_X(k)$, see \cite{CN2} for more precise statements.

In this article, we prove the conjecture in case $d=2$ and $G={\rm GL}_n$ for an arbitrary smooth projective surface:
\begin{thm}[Theorem \ref{403}]
Let $X$ be a smooth projective surface. Let $\shf{M}_X$ be the moduli stack of Higgs bundles of rank $n$ over $X$. Then the image of the Hitchin morphism $h_X: \shf{M}_X \rightarrow \shf{B}_X$ is $\shf{A}_X$, i.e. $sd_X : \shf{M}_X \rightarrow \shf{A}_X$ is surjective.
\end{thm}
The ${\rm GL}_n$ case is the starting point for general results. We believe this approach is applicable to some special reductive groups. However, instead of using the language of ``spectral cover" and working on the problem case by case, we are looking forward to seeing a proof of this conjecture in the language of ``cameral covers" (see \cite{Don1993,DoGai,Falt}, and the idea of the cameral cover in this case is briefly stated in \cite[\S 5]{CN2}).

Here is an overview of the proof of the main result and the structure of the paper. Let ${\rm Chow}^n(T^*_X/X)$ be the relative Chow variety of $n$-points on $T^*_X$ over $X$. Chen and Ng\^o defined the space of spectral data $\mathscr{A}_X$ to be the space of sections $X \rightarrow {\rm Chow}^n(T^*_X/X)$, and a section $a$ is called a \emph{spectral datum}. 

For the relative Chow variety, there is a natural stratification
\begin{align*}
{\rm Chow}^{n}(T^*_X/X) = \coprod_{\mu} {\rm Chow}^n_{\mu}(T^*_X/X),
\end{align*}
where the union is taken over all partitions of $n$. Chen and Ng\^o showed that given a spectral datum $a: X \rightarrow {\rm Chow}^n(T^*_X/X)$, if the generic point of $X$ under $a$ lies in ${\rm Chow}^n_{(1^n)}(T^*_X/X)$, then there is a finite flat cover $\widetilde{X}_a \rightarrow X$ such that the Higgs bundles over $X$ with spectral datum $a$ correspond to the Cohen-Macaulay sheaves of rank one on $\widetilde{X}_a$ (see \cite[Theorems 7.1 and 7.3]{CN2}). In particular, for such general $a$, then $sd_X^{-1}(a)$ is nonempty.

Now given an arbitrary $a \in \mathscr{A}_X(k)$, we assume that $a$ maps the generic point of $X$ to the stratum ${\rm Chow}^n_{\mu}(T^*_X/X)$, where $\mu=(1^{\alpha_1}, 2^{\alpha_2},\dots, s^{\alpha_s})$ is a partition of $n$. Therefore, $a$ induces a morphism $X \rightarrow \overline{{\rm Chow}^n_{\mu}(T^*_X /X)}$. Here $\overline{{\rm Chow}^n_{\mu}(T^*_X /X)}$ denotes the closure of ${\rm Chow}^n_{\mu}(T^*_X/X)$ in ${\rm Chow}^{n}(T^*_X/X)$. Since $X$ is smooth, the spectral datum $a$ factors through the normalization $\overline{{\rm Chow}^n_{\mu}(T^*_X /X)}^{\rm nor}$. Abusing the notation, we use the same notation $a: X \rightarrow \overline{{\rm Chow}^n_{\mu}(T^*_X /X)}^{\rm nor}$. We observe the isomorphism, which is analogous to \cite[Exercise 7.4.2]{FGA},
\begin{align*}
    \overline{{\rm Chow}^n_{\mu}(T^*_X/X)}^{\rm nor} \cong {\rm Chow}^{\alpha_1}(T^*_X/X)\times_X \cdots \times_X {\rm Chow}^{\alpha_s}(T^*_X/X)
\end{align*}
still holds (Proposition \ref{normalization}). It follows that $a$ yields a collection of spectral data $(a_i)_{1\leq i \leq s}$,
\begin{align*}
a_i: X \rightarrow \overline{{\rm Chow}^{\alpha_i}_{(1^{\alpha_i})}(T^*_X/X)}={\rm Chow}^{\alpha_i}(T^*_X/X).
\end{align*}
This gives a decomposition of the given spectral datum $a$, which is the key point to prove the main result.
\begin{thm}[Theorem \ref{decomp}]
Let $X$ a smooth projective surface. Given a spectral datum $a: X \rightarrow {\rm Chow}^{n}(T^*_X /X)$, suppose that the generic point of $X$ is mapped to the stratum ${\rm Chow}^{n}_{\mu}(T^*_X /X)$ with $\mu=(1^{\alpha_1}, 2^{\alpha_2},\dots, s^{\alpha_s})$. Then there exist spectral data $a_i: X \rightarrow \overline{{\rm Chow}^{\alpha_i}_{(1^{\alpha_i})}(T^*_X /X)}$, $i=1, \cdots, s$, such that
\[a=\sum^s_{i=1} ia_i.\]
\end{thm}

Since $a_i$ maps the generic point of $X$ to ${\rm Chow}^{\alpha_i}_{(1^{\alpha_i})}(T^*_X/X)$, there exists a Higgs bundle $(E_i,\phi_i)$ corresponds to $a_i$ according to \cite[Theorem 7.3]{CN2}. Then the direct sum
\begin{align*}
(E,\phi):=\bigoplus_{i=1}^s (E_i,\phi_i)^{\oplus i},
\end{align*}
is a Higgs bundle of rank $n$ corresponding to the spectral datum $a$. This gives the surjectivity of $sd_X: \mathscr{M}_X(k) \rightarrow \mathscr{A}_X(k)$.

In \S 5, we prove the birational invariance of $\shf{A}_X$ for any dimensional smooth variety $X$.
\begin{thm}[Theorem \ref{303}]
Let $\chi: X'\rightarrow X$ be a birational proper morphism between smooth varieties. Then the natural morphism $\shf{A}_X\rightarrow \shf{A}_{X'}$ induced by $\chi$ is an isomorphism as $k$-schemes.
\end{thm}

It is shown in \cite[Proposition 8.1]{CN2} that $\shf{A}_X \cong \shf{B}_X$ for any $\prj{1}$-bundle $X$ over a smooth curve. As a quick application of Theorem \ref{303}, we extend the fact to all ruled surfaces \footnote{The usage of the term ``ruled surface" in \cite{CN2} is different from ours. It seems that in \cite{CN2} it simply means a $\prj{1}$-bundle over a smooth curve, while here it means any smooth surface that is birationally equivalent to a $\prj{1}$-bundle over a smooth curve.}, bypassing the requirement as in \cite{CN2} that each fibre in a fibration is reduced. 

There are a number of interesting questions to be addressed in the future.
\begin{enumerate}
    \item In this article, we have only considered Higgs bundles (${\rm GL}_n$-case). It would be very interesting to generalize the result to $G$-Higgs bundles over higher dimensional varieties in the language of cameral covers. This is of course the full statement of Chen and Ng\^o's conjecture.
    \item The classical Hitchin morphism is given for the moduli space of semistable Higgs bundles, and the conjectural image of the Hitchin morphism is for the moduli stack of Higgs bundles. It is natural to ask if $\shf{A}_X$ is the image of the Hitchin morphism is for the moduli space $\mathcal{M}^{ss}_X(P)$ of Higgs bundles with respect to some Hilbert polynomial $P$?
\end{enumerate}

\vspace{2mm}
\textbf{Acknowledgments}.
The authors would like to thank X. Hu, G. Kydonakis, M. Li, J. Xu and L. Zhao for helpful discussions and conversations. The authors are very grateful to G. Kydonakis for invaluable suggestions on an early draft. We sincerely thank the referee for spotting a gap in the proof of Theorem \ref{decomp}, and for numerous suggestions, which to a large extent improve the exposition of the paper. During the preparation of the paper, the first author was partially supported by the Guangdong Basic and Applied Basic Research Foundation No.~2020A1515010876, and the second author was partially supported by National Key R$\&$D Program of China No.~2022YFA1006600 and by NSFC 12101243.
\vspace{2mm}

\section{Moduli Stacks of Higgs Bundles, Hitchin Base and Spectral data}
In this section, we fix a positive integer $n$ as the rank and discuss three moduli stacks on a smooth quasi-projective variety $X$ over $k$:
\begin{itemize}
\item the moduli stack of Higgs bundles $\shf{M}_X$;
\item the moduli stack of the Hitchin base $\shf{B}_X$;
\item the moduli stack of spectral data $\shf{A}_X$,
\end{itemize}
and they are the main objects we work on. In this paper, we only focus on the Higgs bundles, i.e. ${\rm GL}_n$-Higgs bundles, but some of the constructions can be extended to $G$-Higgs bundles (see \cite{CN,CN2} for details).

Let $\mathfrak{C}^d_{{\rm GL}_n} \subseteq \mathfrak{gl}^d_n$ be the closed subscheme consisting of $d$-tuples $(x_1,\dots,x_d) \in \mathfrak{gl}^d_n$ such that $x_i$ and $x_j$ commutes, i.e. $[x_i,x_j]=0$ for all indices $i,j$. The scheme $\mathfrak{C}^d_{{\rm GL}_n}$ is also known as the commuting scheme. There are two group actions on $\mathfrak{C}^d_{{\rm GL}_n}$:
\begin{itemize}
\item Given $g \in {\rm GL}_n$, the ${\rm GL}_n$-action on $\mathfrak{C}^d_{{\rm GL}_n}$ is defined as
    \begin{align*}
    (x_1,\dots,x_d) \mapsto ({\rm ad}(g)(x_1),\dots,{\rm ad}(g)(x_d)).
    \end{align*}
\item Given $h \in {\rm GL}_d$, the ${\rm GL}_d$-action on $\mathfrak{C}^d_{{\rm GL}_n}$ is defined as
    \begin{align*}
    (x_1,\dots,x_d) \mapsto (x_1,\dots,x_d)h.
    \end{align*}
\end{itemize}

We first consider the ${\rm GL}_n$-action on $\mathfrak{C}^d_{{\rm GL}_n}$. With the same argument of the Chevalley restriction theorem, we obtain a ${\rm GL}_d$-equivariant morphism
\begin{align*}
k[\mathfrak{gl}_n^d]^{{\rm GL}_n} \rightarrow k[\mathfrak{t}^d]^{\mathfrak{S}_n},
\end{align*}
following the inclusion $\mathfrak{t}^d \rightarrow \mathfrak{gl}_n^d$, where $\mathfrak{t} \subseteq \mathfrak{gl}_n$ is a Cartan subalgebra and $\mathfrak{S}_n$ is the permutation group. Restricting the morphism to $\mathfrak{C}^d_{{\rm GL}_n}$, we have a ${\rm GL}_d$-equivariant lifting
\begin{align*}
[\mathfrak{C}^{d}_{{\rm GL}_n}/{\rm GL}_n] \rightarrow \Spec(k[\mathfrak{t}^d]^{\mathfrak{S}_n})
\end{align*}
of $[\mathfrak{C}^{d}_{{\rm GL}_n}/{\rm GL}_n] \rightarrow \mathfrak{C}^{d}_{{\rm GL}_n}/{\rm GL}_n$. Let $V$ be a vector space of dimension $d$. We have
\begin{align*}
\Spec(k[\mathfrak{t}^d]^{\mathfrak{S}_n}) \cong {\rm Chow}^n(V),
\end{align*}
where ${\rm Chow}^n(V)$ is the Chow variety of $n$ points on $V$. Then we have the morphism
\begin{align*}
[\mathfrak{C}^{d}_{{\rm GL}_n}/{\rm GL}_n] \rightarrow {\rm Chow}^n(V).
\end{align*}
By \cite[Lemma 2.1]{CN}, we know that there is a closed embedding
\begin{align*}
{\rm Chow}^n(V) & \hookrightarrow V \times S^2 V \times \dots \times S^n V, \\
[x_1,\dots,x_n] & \mapsto (c_1,\dots,c_n),
\end{align*}
where $c_i \in S^i V$ is the $i$-th elementary symmetric polynomial in variables $x_1,\dots,x_n \in V$. This induces the following morphism
\begin{align*}
[\mathfrak{C}^{d}_{{\rm GL}_n}/{\rm GL}_n] \rightarrow {\rm Chow}^n(V)\hookrightarrow V \times S^2 V \times \dots \times S^n V.
\end{align*}

As a $d$-dimensional vector space, there is a canonical ${\rm GL}_d$-action on $V$. We define the quotient stacks
\begin{align*}
\shf{A}:=[{\rm Chow}^n(V)/{\rm GL}_d], \quad \shf{B}:= [V \times S^2V \times \dots \times S^n V/ {\rm GL}_d].
\end{align*}
By the definition of quotient stacks, for each $k$-scheme $S$, $\shf{B}(S)$ is the groupoid of pairs $(E,\phi)$, where $E$ is a rank $d$ vector bundle (as a ${\rm GL}_d$-principal bundle) over $S$ and
\begin{align*}
\phi: E \rightarrow V \times S^2V \times \dots \times S^n V
\end{align*}
is a ${\rm GL}_d$-equivariant map. Note that the ${\rm GL}_d$-equivariant map $\phi$ corresponds to a unique element in $\bigoplus\limits_{i=1}^{n}H^0(S, S^i E)$. Therefore, we have the following isomorphism as sets (or groupoids) by restricting to a fixed vector bundle $E$ over $S$,
\begin{align*}
\shf{B}(S)|_E \cong \bigoplus_{i=1}^{n}H^0(S, S^i E).
\end{align*}
This isomorphism also implies that $\shf{B}(S)|_E$ has a scheme structure.

For the quotient stack $\shf{A}$, given a $k$-scheme $S$, $\shf{A}(S)$ is the groupoid of pairs $(E,\varphi)$, where $E$ is a rank $d$ vector bundle over $S$ and
\begin{align*}
\varphi: E \rightarrow {\rm Chow}^n(V)
\end{align*}
is a ${\rm GL}_d$-equivariant map. This map is equivalent to an $S$-point in ${\rm Chow}^n(E/S)$, where ${\rm Chow}^n(E/S)$ is the relative Chow variety. In other words, we have
\begin{align*}
\shf{A}(S)|_E \cong {\rm Sect}(S,{\rm Chow}^n(E/S)),
\end{align*}
where ${\rm Sect}(S,{\rm Chow}^n(E/S))$ is the set of sections.

\subsection*{Moduli stack of Higgs bundles}
Recall that there is a natural ${\rm GL}_d$-action on $[\mathfrak{C}^{d}_{{\rm GL}_n}/{\rm GL}_n]$. Denote by $\shf{M}:=[\mathfrak{C}^{d}_{{\rm GL}_n}/{\rm GL}_n \times {\rm GL}_d]$ the quotient stack. Let $S$ be a scheme, and then the $S$-points of $\shf{M}$ are triples $(V,E,\theta)$, where
\begin{itemize}
\item $V$ is a vector bundle of rank $d$ (${\rm GL}_d$-principal bundle) over $S$;
\item $E$ is a vector bundle of rank $n$ (${\rm GL}_n$-principal bundle) over $S$;
\item $\theta \in H^0(S, \mathcal{E}nd(E) \otimes V)$ is a $\mathcal{O}_S$-linear morphism such that $\theta \wedge \theta=0$.
\end{itemize}

Now we fix a smooth quasi-projective variety $X$ of dimension $d$, and let $T^*_X$ denote the cotangent bundle over $X$. Let 
\begin{align*}
    \sigma: X \rightarrow [\ast / {\rm GL}_d] = \mathbb{B}{\rm GL}_d
\end{align*}
be the morphism corresponding to $T^*_X$. Also, there is a natural morphism
\begin{align*}
\shf{M}\rightarrow \mathbb{B}{\rm GL}_d.
\end{align*}
We consider the fiber product 
\begin{center}
\begin{tikzcd}
    \shf{M} \times_{\mathbb{B}{\rm GL}_d} X \arrow[r] \arrow[d] &\shf{M} \arrow[d] \\
    X  \arrow[r,"\sigma"] & \mathbb{B}{\rm GL}_d.
\end{tikzcd}
\end{center}
Each section $X \rightarrow \shf{M} \times_{\mathbb{B}{\rm GL}_d} X$ gives a Higgs bundle $(E,\theta)$ on $X$, where
\begin{itemize}
\item $E$ is a rank $n$ vector bundle over $X$;
\item $\theta \in H^0(X, \mathcal{E}nd(E) \otimes T^*_X)$ such that $\theta \wedge \theta=0$.
\end{itemize}
Denote by 
\begin{align*}
    \shf{M}_X:={\rm Sect}(X, \shf{M} \times_{\mathbb{B}{\rm GL}_d} X),
\end{align*}
the \emph{moduli stack of Higgs bundles}.

\subsection*{Hitchin base $\shf{B}_X$}
We have a natural morphism $\shf{B} \rightarrow \mathbb{B}{\rm GL}_d$ of stacks. Similarly, we define
\begin{align*}
    \shf{B}_X:= {\rm Sect}(X, \shf{B} \times_{\mathbb{B}{\rm GL}_d} X)
\end{align*}
to be the stack of sections of the fiber product $\shf{B} \times_X \mathbb{B}{\rm GL}_d$. The stack $\shf{B}_X$ is actually a space:
\begin{lem}\label{201}
The stack $\shf{B}_X$ is represented by $\bigoplus\limits_{i=1}^{n} H^0(X, S^i T^*_X)$.
\end{lem}
Abusing the notation, we also use $\shf{B}_X$ for  $\bigoplus\limits_{i=1}^{n} H^0(X, S^i T^*_X)$, which is called the \emph{Hitchin base}.

\subsection*{Space of spectral data $\shf{A}_X$}
Similar to $\shf{B}_X$, we define $\shf{A}_X$ as
\begin{align*}
    \shf{A}_X:= {\rm Sect}(X, \shf{A} \times_{\mathbb{B}{\rm GL}_d} X)
\end{align*}
We have the following lemma.
\begin{lem}\label{202}
The stack $\shf{A}_X$ is represented by the scheme ${\rm Sect}(X,  {\rm Chow}^n( T^*_X/X ))$. 
\end{lem}
We use the same notation $\shf{A}_X$ for the corresponding scheme. An element $a \in \shf{A}_X$ is called \emph{a spectral datum} and $\shf{A}_X$ is called \emph{the space of spectral data}. Sometimes, we use the notation $\shf{A}_X^n$ to highlight the rank $n$.

\begin{rem}\label{204}
Let $U \rightarrow X$ be a morphism of schemes. Since $\shf{M}$ is a stack, we have a natural map $\shf{M}(X) \rightarrow \shf{M}(U)$. Furthermore, this map induces a morphism $\shf{M}_X \rightarrow \shf{M}_U$ as stacks. The same argument also works for morphisms $\shf{A}_X \rightarrow \shf{A}_U$ and $\shf{B}_X \rightarrow \shf{B}_U$. As a special case, if $U$ is an open subscheme of $X$, then we have a natural morphism $\shf{A}_X \rightarrow \shf{A}_U$.
\end{rem}

\begin{rem}\label{205}
Recall that the closed embedding
\begin{align*}
{\rm Chow}^n(V) \hookrightarrow V \times S^2 V \times \dots \times S^n V,
\end{align*}
induces the closed embedding $\mathfrak{i}_X: \shf{A}_X \hookrightarrow \shf{B}_X$ (\cite[\S 6 and \S 7]{CN2}). Chen and Ng\^o showed that the Hitchin morphism $h_X: \shf{M}_X \rightarrow \shf{B}_X$ factors through $\shf{A}_X$. More precisely, there exists a map $sd_X: \shf{M}_X \rightarrow \shf{A}_X$ such that $h_X=\mathfrak{i}_X \circ sd_X$.
\begin{center}
\begin{tikzcd}
            & \shf{M}_X \arrow[d,"h_X"] \arrow[ld, dotted, "sd_X" description]  \\
            \shf{A}_X \arrow[r, hook, "\mathfrak{i}_X"] & \shf{B}_X
\end{tikzcd}
\end{center}
The morphism $sd_X$ is called \emph{the spectral data morphism}.
\end{rem}

\section{Decomposition of spectral data}

In this section, we will prove a decomposition theorem (Theorem \ref{decomp}) for spectral datum $a \in \shf{A}_X$, where $X$ is an algebraic surface. This decomposition theorem is an important tool to prove the main result Theorem \ref{403}. In \S 3.1, we briefly review some relevant facts about Chow varieties, and refer the reader to \cite[Chap. 7]{FGA} for more details. In \S 3.2, we generalize the results to the relative case, and prove the decomposition theorem.

\subsection{Classical case}
Let $V$ be a quasi-projective surface and let ${\rm Chow}^n(V)$ be the Chow variety of $n$ points. For a partition $\mu$ of $n$, we will write $\mu=(n_1,\dots,n_k)$ or $\mu=(1^{\alpha_1}, 2^{\alpha_2},\dots, s^{\alpha_s})$ interchangeably. Define a locally closed subset ${\rm Chow}^n_{\mu}(V)$ of ${\rm Chow}^n(V)$ as
\begin{align*}
{\rm Chow}^n_{\mu}(V)=\{\sum_i^k n_i z_i \: |\:  z_i \in V \text{ distinct points }  \}.
\end{align*}
Then there is a stratification:
\begin{align*}
{\rm Chow}^{n}(V) = \coprod_{\mu} {\rm Chow}^n_{\mu}(V),
\end{align*}
where $\mu$ ranges over all partitions of $n$.

Given a partition $\mu=(1^{\alpha_1}2^{\alpha_2}\dots s^{\alpha_s})$, the natural morphism
\begin{align*}
{\rm Chow}^{\alpha_1}(V) \times {\rm Chow}^{\alpha_2}(V) \times \dots \times {\rm Chow}^{\alpha_s}(V) \rightarrow {\rm Chow}^n(V)
\end{align*}
defined by
\begin{align*}
(z_1,z_2,\dots , z_s) \mapsto z_1+ 2z_2 + \dots + s z_s.
\end{align*}
indeed gives the normalization of $\overline{{\rm Chow}^n_{\mu}(V)}$.
\begin{lem}[Exercise 7.4.2 in \cite{FGA}]\label{101}
Given a partition $\mu=(1^{\alpha_1}, 2^{\alpha_2}, \dots, s^{\alpha_s})$, we have
\begin{align*}
\overline{{\rm Chow}^n_{\mu}(V)}^{\rm nor} \cong {\rm Chow}^{\alpha_1}(V) \times \dots \times {\rm Chow}^{\alpha_s}(V),
\end{align*}
where $\overline{{\rm Chow}^n_{\mu}(V)}^{\rm nor}$ is the normalization of the closure of ${\rm Chow}^n_{\mu}(V)$.
\end{lem}

\subsection{Relative case}
Let $X$ be a smooth projective surface. Denote by $T^*_X$ the cotangent bundle over $X$. The Chow variety ${\rm Chow}^n(T^*_X/X)$ is defined as follows
\begin{align*}
{\rm Chow}^n(T^*_X/X):= \underbrace{T^*_X \times_X \dots \times_X T^*_X}_n/ \mathfrak{S}_n,
\end{align*}
where $\mathfrak{S}_n$ is the symmetric group of $n$ letters. Similar to the classical case, one can define a locally closed subset ${\rm Chow}^n_{\mu}(T^*_X/X)$ of ${\rm Chow}^n(T^*_X/X)$ with respect to the partition $\mu$, and the stratification of ${\rm Chow}^{n}(T^*_X/X)$ holds
\begin{align*}
{\rm Chow}^{n}(T^*_X/X) = \coprod_{\mu} {\rm Chow}^n_{\mu}(T^*_X/X),
\end{align*}
where $\mu$ ranges over all partitions of $n$.

For ease of notation, set
\begin{align*}
\prod^s_{j=1, X}{\rm Chow}^{\alpha_j}(T^*_X/X):={\rm Chow}^{\alpha_1}(T^*_X/X) \times_X {\rm Chow}^{\alpha_2}(T^*_X/X) \times_X \dots \times_X {\rm Chow}^{\alpha_s}(T^*_X/X),
\end{align*}
where the subscript $X$ under the symbol $\prod^s_{j=1, X}$ indicates that the product is taken over $X$. Similarly set
\begin{align*}
\prod^n_{j=1, X}T^*_X:=\underbrace{T^*_X \times_X \dots \times_X T^*_X}_n.
\end{align*}

We summarize some easy facts:
\begin{lem}\label{102}
For any $n\ge 1$ and a partition $\mu=(1^{\alpha_1}, 2^{\alpha_2},\dots, s^{\alpha_s})$, the following hold:
\begin{itemize}
    \item[(i)] the scheme $\prod\limits^n_{j=1, X}T^*_X$ is smooth over $X$;
    \item[(ii)] the scheme ${\rm Chow}^n_{\mu}(T^*_X /X)$ is smooth over $X$;
    \item[(iii)] the scheme ${\rm Chow}^n(T^*_X /X)$ has rational singularities, in particular it is normal and Cohen-Macaulay;
    \item[(iv)] the product $\prod\limits^s_{j=1, X}{\rm Chow}^{\alpha_j}(T^*_X/X)\cong \prod\limits^{\alpha_1+\dots+\alpha_s}_{j=1, X}T^*_X/{\mathfrak{S}_{\alpha_1}\times\cdots \times\mathfrak{S}_{\alpha_s}}$ has rational singularities, in particular, it is normal and Cohen-Macaulay.
\end{itemize}
\end{lem}
\begin{proof}
(i) follows from the facts that $T^*_X$ is smooth over $X$ and smoothness is stable under base change. (ii) Take a trivialization of $T^*_X$ over an open covering $\{U_i\}$ of $X$. Over each $U_i$, ${\rm Chow}^n_{\mu}(T^*_X /X)$ is isomorphic to ${\rm Chow}^n_{\mu}(\mathbb{A}^2)\times_k U_i$. Since the stratum ${\rm Chow}^n_{\mu}(\mathbb{A}^2)$ is smooth over $k$, so is ${\rm Chow}^n_{\mu}(T^*_X /X)$ smooth over $U_i$. (iii) is true because ${\rm Chow}^n(T^*_X /X)$ is the quotient of the smooth variety $\prod\limits^n_{j=1, X}T^*_X$ by a finite group. (iv) is analogues to (iii).
\end{proof}

Given a partition $\mu=(1^{\alpha_1}, 2^{\alpha_2},\dots, s^{\alpha_s})$ of $n$, we define a morphism
\begin{equation*}
\tau_{\mu}: \prod^s_{j=1, X}{\rm Chow}^{\alpha_j}(T^*_X/X) \rightarrow {\rm Chow}^n(T^*_X/X)
\end{equation*}
as
\begin{align*}
(z_1,z_2,\dots , z_s) \mapsto  z_1+ 2z_2 + \dots + s z_s.
\end{align*}
The following is a key observation, which is a relative version of Lemma \ref{101}.

\begin{prop}\label{normalization}
Given a partition $\mu=(1^{\alpha_1}, 2^{\alpha_2},\dots, s^{\alpha_s})$, the morphism $\tau_{\mu}$ gives the normalization of $\overline{{\rm Chow}^n_{\mu}(T^*_X/X)}$.
\end{prop}

\begin{proof}
The generic point of the product is mapped to the generic point of ${\rm Chow}^n_{\mu}(T^*_X/X)$ by $\tau_{\mu}$. Thus $\tau_{\mu}$ induces a proper, surjective morphism from $\prod^s_{j=1, X}{\rm Chow}^{\alpha_j}(T^*_X/X)$ to
$\overline{{\rm Chow}^n_{\mu}(T^*_X/X)}$. Also, note that by Lemma \ref{102} (iv), the product $\prod^s_{j=1, X}{\rm Chow}^{\alpha_j}(T^*_X/X)$
is a normal variety.

To see $\tau_{\mu}$ is indeed the normalization, we shall utilize Lemma \ref{103}, which says one can check normalization locally. Take an open covering $\{U_i\}$ of $X$ such that $T^*_X |_{U_i} \cong U_i \times \mathbb{A}^2$. We have the commutative diagram with Cartesian squares
$$\xymatrix{
&\prod^s_{j=1, U_i}{\rm Chow}^{\alpha_j}(T^*_{U_i}/{U_i})\ar[d]^{\tau'_{\mu}}\ar[ld]\ar[rrd]&&\\
\prod^s_{j=1, X}{\rm Chow}^{\alpha_j}(T^*_X/X)\ar[rrd]\ar[d]^{\tau_{\mu}}& \overline{{\rm Chow}^n_{\mu}(T^*_{U_i}/{U_i})}\ar@{^(->}[rr] \ar[ld]\ar[rrd] &  & {\rm Chow}^n(T^*_{U_i}/{U_i})\ar[d]\ar[ld]\\
\overline{{\rm Chow}^n_{\mu}(T^*_X/X)}\ar@{^(->}[rr]\ar[rrd]&  &  {\rm Chow}^n(T^*_X/X)\ar[d] & U_i\ar[ld], \\
& & X &
}$$
where $\overline{{\rm Chow}_{\mu}^n(T^*_{U_i} /U_i)}$ is the closure of ${\rm Chow}_{\mu}^n(T^*_{U_i} /U_i)$ in ${\rm Chow}^n(T^*_{U_i}/{U_i})$. We have the isomorphisms
\begin{align*}
\overline{{\rm Chow}_{\mu}^n(T^*_{U_i} /U_i)}  & \cong  \overline{{\rm Chow}_{\mu}^n(\mathbb{A}^2)} \times U_i,\\
\prod^s_{j=1, U_i}{\rm Chow}^{\alpha_j}(T^*_{U_i}/U_i) & \cong \big(\prod^s_{j=1}{\rm Chow}^{\alpha_j}(\mathbb{A}^2)\big)\times U_i.
\end{align*}
Thus by Lemma \ref{101}, $\tau'_{\mu}$ is the normalization of $\overline{{\rm Chow}^n_{\mu}(T^*_{U_i}/{U_i})}$. Applying Lemma \ref{103}, we get $\tau_{\mu}$ is the normalization of $\overline{{\rm Chow}^n_{\mu}(T^*_X/X)}$.
\end{proof}

\begin{lem}\label{103}
Let $X$ be an integral scheme and $\tilde{X}$ a normal integral scheme. Given a morphism $\nu: \tilde{X}\rightarrow X$, suppose there exists an open covering $\{U_i\}$ of $X$ such that $\nu_i: \tilde{U_i}=\nu^{-1}(U_i)\rightarrow U_i$ is the normalization for each $i$. Then $\nu$ is the normalization.
\end{lem}
\begin{proof}
This follows from the universal property of normalization \cite[II Ex.~3.8]{Har}.
\end{proof}

The following strengthened universal property of normalization should be well-known to the experts. It is stated in a more general form than what we actually need. For the convenience of the reader, we include a proof here, which follows closely an argument of Johann Haas \cite{Haas}. 
\begin{lem}\label{104}
Let $f: Y\rightarrow X$ be a morphism between integral schemes. Let $\nu :X^{\nu}\rightarrow X$ be the normalization of $X$.
Assume $X$ is excellent. Suppose $Y$ is normal and $f$ maps the generic point of $Y$ to a normal point of $X$. Then there exists a unique morphism $f': Y\rightarrow X^{\nu}$ such that $\nu\circ f'=f$.
\end{lem}

\begin{proof}
Let $Z$ be the scheme-theoretic image of $f$. As a closed subscheme of $X$, $Z$ is excellent. Since $f: Y\rightarrow Z$ is dominant, $f$ factors uniquely through the normalization $Z^{\nu}$ of $Z$ \cite[II Ex. 3.8]{Har}. So after replacing $Y$ by $Z^{\nu}$, we can assume that $f$ induces the normalization of its image $Z$. 

Consider the diagram 
\begin{equation}\label{norm}
\xymatrix{
 \overline{Z}\ar@{^{(}->}[r] \ar[rd]^{\mu} & X^{\nu}\times_X Z\ar@{^{(}->}[r]\ar[d] & X^{\nu}\ar[d]^{\nu}\\
Y\ar[r]^{f} & Z\ar@{^{(}->}[r] & X}
\end{equation}
with the Cartesian square. Since the generic point of $Z$ is a normal point of $X$, and since the locus of normal points of $X$ is open \cite[7.8.3]{Gro} due to the excellence of $X$, there exists a unique irreducible component $\overline{Z}$ of $X^{\nu}\times_X Z$ dominating $Z$. Let $\overline{Z}$ be given the reduced subscheme structure. Then $\overline{Z}$ is an integral scheme and the induced $\mu: \overline{Z}\rightarrow Z$ is finite and birational.  

All morphisms involved are affine, so we can assume $Z=\Spec A$, $\overline{Z}=\Spec{B}$, and $Y=\Spec{\tilde{A}}$, where $\tilde{A}$ is the integral closure of $A$ in the field $K(A)$ of fractions. Moreover $B$ is integral over $A$, and the induced map $K(A)\rightarrow K(B)$ by $\mu$ is an isomorphism of fields. Clearly there exists a unique ring homomorphism $h: B\rightarrow \tilde{A}$ fitting into the diagram
\begin{equation}\label{unique}
\xymatrix{
K(A)\ar[r] & K(B)\\
\tilde{A}\ar@{^{(}->}[u] & \\
A\ar@{^{(}->}[u]\ar@{^{(}->}[r] & B\ar@{^{(}->}[uu]\ar@{-->}[lu]_{h}}
\end{equation}
This gives a morphism $Y\rightarrow \overline{Z}$ fitting into  diagram (\ref{norm}) and hence $f': Y\rightarrow X^{\nu}$. 

To see the uniqueness, note $f'$ with the property $\nu\circ f'=f$ will factor through $X^{\nu}\times_X Z$ uniquely. As $f$ has been assumed to be dominant after reduction, $f'$ further factors through $\overline{Z}$ uniquely. Then the uniqueness of $f'$ results from that of $h$ in diagram \ref{unique}. 
\end{proof}

Now we are ready to prove the decomposition theorem for spectral data.
\begin{thm}\label{decomp}
Let $X$ be a smooth projective surface. Given a spectral datum $a: X \rightarrow {\rm Chow}^{n}(T^*_X /X)$, suppose that the generic point of $X$ is mapped to the stratum ${\rm Chow}^{n}_{\mu}(T^*_X /X)$ with $\mu=(1^{\alpha_1}, 2^{\alpha_2},\dots, s^{\alpha_s})$. Then there exist spectral data $a_i: X \rightarrow \overline{{\rm Chow}^{\alpha_i}_{(1^{\alpha_i})}(T^*_X /X)}$, $i=1, \cdots, s$, such that
\[a=\sum^s_{i=1} ia_i.\]
\end{thm}
\begin{proof}
By the assumption, $a$ induces a morphism $\bar{a}: X \rightarrow \overline{{\rm Chow}^n_{\mu}(T^*_X /X)}$ over $X$. Note $\overline{{\rm Chow}^n_{\mu}(T^*_X /X)}$ is an algebraic variety, in particular locally of finite type over a field $k$, so is excellent \cite[7.8.3]{Gro}.
Since $\bar{a}$ maps the generic point of $X$ to a smooth point of $\overline{{\rm Chow}^n_{\mu}(T^*_X /X)}$, see Lemma \ref{102}(ii), it follows from Lemma \ref{104} that
there exists a unique $a'$ fitting into the diagram below
\begin{center}
\begin{tikzcd}
    & \overline{{\rm Chow}^n_{\mu}(T^*_X /X)}^{\rm nor} \arrow[d] \arrow[rd, "\tau_{\mu}"]  & \\
    X \arrow[r, "\bar{a}"]\arrow[ru,"a'", dotted]  & \overline{{\rm Chow}^n_{\mu}(T^*_X /X)} \arrow[r, hook] & {\rm Chow}^{n}(T^*_X /X),
\end{tikzcd}
\end{center}
where by Proposition \ref{normalization}, 
\begin{align*}
\overline{{\rm Chow}^n_{\mu}(T^*_X /X)}^{\rm nor} \cong \prod^s_{j=1, X}{\rm Chow}^{\alpha_j}(T^*_X/X).
\end{align*}
It turns out $a'$ is equivalent to data $(a_i)_{1 \leq i \leq s}$, where
\begin{align*}
a_i: X \rightarrow  \overline{{\rm Chow}^{\alpha_i}_{(1^{\alpha_i})}(T^*_X /X)}.
\end{align*}
Finally by the commutativity of the diagram, we conclude that $a=\sum\limits_{i=1}^s i a_i$.
\end{proof}

\section{The image of the Hitchin morphism for surfaces}

Given a spectral datum $a: X \rightarrow {\rm Chow}^{n}(T^*_X/X)$, we first briefly review the construction of the spectral cover $\widetilde{X}_a \rightarrow X$, and the correspondence between Cohen-Macaulay sheaves of rank one over $\widetilde{X}_a$ and Higgs bundles over $X$ \cite[\S 7]{CN2}. Next, we show that given two Higgs bundles, the direct sum corresponds to the sum of the corresponding spectral data (Lemma \ref{402}). Finally, we prove the main theorem that $sd_X$ is surjective (Theorem \ref{403}) with the help of Theorem \ref{decomp}.

We first review the construction of Cayley scheme ${\rm Cayley}^n(V)$ (see \cite[\S 6]{CN2} for instance). Consider the morphism
\begin{align*}
\chi : {\rm Chow}^n(V)  \times V \rightarrow S^n V   
\end{align*}
defined as
\begin{align*}
    ([x_1,\dots,x_n],x) \mapsto x^n - c_1 x^{n_1} + \dots + (-1)^n c_n,
\end{align*}
where $c_i$ is the $i$-th elementary symmetric polynomial of variables $x_1,\dots,x_n$. The closed subcheme ${\rm Cayley}^n(V)$ is defined as $\chi^{-1}(0)$. As a relative version, ${\rm Cayley}^n(T^*_X/X)$ is a closed subscheme of ${\rm Chow}^n(T^*_X/X) \times_X T^*_X$. Then $X_a$ is defined as the pullback in the following diagram
\begin{center}
\begin{tikzcd}
    X_a \arrow[d,"\pi_a"] \arrow[rr]& & {\rm Cayley}^n(T^*_X/X) \arrow[d,"\pi"] \\
    X  \arrow[rr,"a"] & & {\rm Chow}^{n}(T^*_X/X).
\end{tikzcd}
\end{center}

\begin{thm}[Theorems 7.1 and 7.3 in \cite{CN2}]\label{401}
Given a spectral datum $a: X \rightarrow {\rm Chow}^{n}(T^*_X/X)$ such that $a$ maps the generic point of $X$ into ${\rm Chow}^n_{(1^n)}(T^*_X/X)$, there exists a unique finite, flat covering $\widetilde{\pi}_a: \widetilde{X}_a \rightarrow X$ such that
\begin{enumerate}
\item there exists an open subset $U \subseteq X$ of codimension at least $2$ such that for every point $x \in U$, the fiber $(\widetilde{\pi}_a)^{-1}(x)$ is a point of ${\rm Hilb}^n(T^*_X/X)$ lying over the point $a(x) \in {\rm Chow}^n(T^*_X/X)$;
\item there is a natural morphism $\widetilde{\iota}_a: \widetilde{X}_a \rightarrow T^*_X$ factoring through $\iota_a: X_a \hookrightarrow T^*_X$ such that the following diagram is commutative;
    \begin{center}
    \begin{tikzcd}
        \widetilde{X}_a \arrow[rd,"\widetilde{\pi}_a"] \arrow[r, "\iota"] & X_a \arrow[r,"\iota_a", hook] \arrow[d, "\pi_a"] & T^*_X \arrow[ld, "\pi"]\\
        & X &
    \end{tikzcd}
    \end{center}
\item $\widetilde{X}_a$ is the Cohen-Macaulayfication of $X_a$;
\item the fiber $sd_X^{-1}(a)$ is isomorphic to the stack of Cohen-Macaulay sheaves of generic rank one over $\widetilde{X}_a$.
\end{enumerate}
\end{thm}
The correspondence demonstrated in the theorem can be understood in the following way. Given a Cohen-Macaulay sheaf $L$ of generic rank one over $\widetilde{X}_a$, its pushfoward $\iota_* L$ is a sheaf over $X_a$. Via $\iota_a$, the sheaf $(\widetilde{\iota}_a)_* L$ has a natural structure as a $S(T_X)$-module, where $S(T_X)$ is the symmetric product of the tangent bundle $T_X$. Furthermore, it is a finite $S(T_X)$-module. Therefore, it corresponds to a spectral data in $\shf{A}_X$, which is exactly $a$, and also corresponds to a Higgs bundle $(E,\phi)$ over $X$.

The correspondence between a finite $S(T_X)$-module $L$ and a spectral data $a$ can be understood in the following way. For each $x \in X$, we can equip $L|_{\pi^{-1}(x)}$ with a cycle
\begin{align*}
{\rm cyc}(L|_{\pi^{-1}(x)}):=\sum_{y \in \pi^{-1}(x)} {\rm len}_{y}(L|_{\pi^{-1}(x)}) \cdot y.
\end{align*}
Then, the $S(T_X)$-module $L$ corresponds to $a$ if for each $x \in X$, we have
\begin{align*}
a(x)={\rm cyc}(L|_{\pi^{-1}(x)}).
\end{align*}

Theorem \ref{401} only works for spectral data $a \in \shf{A}_X$ which map the generic point of $X$ into ${\rm Chow}^n_{(1^n)}(T^*_X/X)$. Given an arbitrary spectral datum, we can follow the approach in the theorem to construct $X_a$, but the correspondence is not clear. We will discuss this issue in Remark \ref{404}.

\begin{lem}\label{402}
Given a pair of spectral data $a_i: X\rightarrow {\rm Chow}^{n_i}(T^*_X/X)$, $i=1, 2$, let $(E_i, \phi_i)$ be Higgs bundles on $X$ whose spectral data are $a_i$. Then $(E_1\oplus E_2, \phi_1\oplus \phi_2)$ is a Higgs bundle with the spectral datum $a_1+a_2: X\rightarrow {\rm Chow}^{n_1+n_2}(T^*_X/X)$.
\end{lem}

\begin{proof}
The pair $(E_1\oplus E_2, \phi_1\oplus \phi_2)$ is obviously a Higgs bundle on $X$. By the correspondence \cite[Lemma 6.8]{Simp2}, there exist coherent sheaves $F_i$ on $T^*_X$ such that $\pi_* F_i\iso E_i$. Thus, $\pi_*(F_1\oplus F_2)\iso E_1\oplus E_2$. For any $x\in X$,
\begin{align*}
a_1(x)+a_2(x)={\rm cyc}( L_1 |_{\pi^{-1}(x)})+{\rm cyc}(L_2|_{\pi^{-1}(x)})={\rm cyc}((L_1+L_2)|_{\pi^{-1}(x)})=(a_1+a_2)(x).
\end{align*}
This finishes the proof.
\end{proof}

\begin{thm}\label{403}
Let $X$ be a smooth projective surface. The image of the Hitchin map $h_X: \shf{M}_X \rightarrow \shf{B}_X$ is $\shf{A}_X$, i.e. $sd_X : \shf{M}_X \rightarrow \shf{A}_X$ is surjective.
\end{thm}

\begin{proof}
We will show that given any spectral datum $a \in \shf{A}_X$, we can construct a Higgs bundle $(E,\phi) \in sd_X^{-1}(a) \subseteq \shf{M}_X$. By Theorem \ref{401}, this is equivalent to finding a $S(T_X)$-module $L$ corresponding to the spectral datum $a$. Now given a spectral datum $a: X \rightarrow {\rm Chow}^n(T^*_X/X)$, suppose that the generic point of $X$ is mapped into some stratum ${\rm Chow}^n_{\mu}(T^*_X/X)$, where $\mu=(1^{\alpha_1}, 2^{\alpha_2}, \dots, s^{\alpha_s})$. By Theorem \ref{decomp}, we have
\begin{align*}
a=\sum^s_{i=1} ia_i,
\end{align*}
where the spectral datum $a_i$ maps the generic point of $X$ into ${\rm Chow}^{\alpha_i}_{(1^{\alpha_i})}(T^*_X /X)$. Therefore, by Theorem \ref{401}, for each $i$, one can construct a spectral cover $\widetilde{X}_{a_i}$ of $X$
\begin{center}
    \begin{tikzcd}
        \widetilde{X}_{a_i} \arrow[rd,"\widetilde{\pi}_i"] \arrow[r] & X_{a_i} \arrow[r,"\iota_i", hook] \arrow[d, "\pi_i"] & T^*_X \arrow[ld, "\pi"]\\
        & X &
    \end{tikzcd}
    \end{center}
such that $\widetilde{\pi}_{a_i}: \widetilde{X}_{a_i} \rightarrow X$ is a finite flat morphism and the morphism $\widetilde{\iota}_{a_i}: \widetilde{X}_{a_i} \rightarrow T^*_X$ factors through $\iota_{i}: X_{a_i} \hookrightarrow T^*_X$. Then, the $S(T_X)$-module $(\widetilde{\iota}_i)_* \mathcal{O}_{\widetilde{X}_{a_i}}$ corresponds to the spectral datum $a_i$. By Lemma \ref{402}, the sheaf
\begin{align*}
L:=\bigoplus\limits_{i=1}^s (\widetilde{\iota}_i)_* \mathcal{O}_{\widetilde{X}_{a_i}}^{\oplus i}
\end{align*}
corresponds to the spectral data $\sum\limits_{i=1}^s i a_i$, which is exactly $a$. Therefore, $sd_X^{-1}(a)$ is nonempty.
\end{proof}

\begin{rem}\label{404}
Although Theorem \ref{403} gives the surjectivity of the morphism $sd_X: \mathscr{M}_X \rightarrow \mathscr{A}_X$, we do not find a cover $\widetilde{X}_a \rightarrow X$ such that Higgs bundles over $X$ with spectral datum $a$ will correspond to some special sheaves over $\widetilde{X}_a$. 
\end{rem}

\section{Birational invariance of the space of spectral data}
The purpose of this section is to show that for a smooth variety $X$ of any dimension, the space of spectral data $\shf{A}_X$ is a birational invariant. Then we apply the result to show the space of spectral data for \textit{all} ruled surfaces are affine spaces. 

\begin{lem}\label{301}
Let $\chi: X'\rightarrow X$ be a birational morphism of smooth varieties. Then the natural morphism $\shf{A}_X\rightarrow \shf{A}_{X'}$ is a closed immersion of schemes.
\end{lem}
\begin{proof}
We have the commutative diagram
$$\xymatrix{
\shf{A}_X \ar[r]\ar@{^(->}[d]^{i_X}&  \shf{A}_{X'}\ar@{^(->}[d]^{i_{X'}}  \\
\shf{B}_X\ar[r]^{\chi^*}&  \shf{B}_{X'},
}$$
where $\mathfrak{i}_X$ and $\mathfrak{i}_{X'}$ are closed immersions by \cite[Lemma 2.1]{CN} and
\[\chi^*: \bigoplus^n_{i=1} H^0(X, S^i T^*_X)\rightarrow \bigoplus^n_{i=1} H^0(X', S^i T^*_{X'})\]
is an injection. From the diagram, we conclude that the natural morphism is a closed immersion.
\end{proof}

\begin{lem}\label{302}
Let $U$ be an open subset of a smooth variety $X$ with $\text{codim}(X\backslash U, X)\ge 2$. Then $\shf{A}_X\rightarrow \shf{A}_U$ is an isomorphism.
\end{lem}

\begin{proof}
Recall for any $k$-scheme $T$, the $T$-points of $\shf{A}_X$ are
\begin{equation*}
    \shf{A}_X(T)=\big\{\xymatrix{
X_T\ar[r] \ar[rd]_{id_{X_T}} &  \text{Chow}^n(T^*_{X_T/T}/{X_T})\ar[d] \\
& X_T,
}\big\}
\end{equation*}
By Lemma \ref{301} and the Yoneda Lemma, we shall show that the natural map $\shf{A}_X(T)\rightarrow \shf{A}_U(T)$ induced by $j$ is surjective for all such $T$.

First we consider the case $T^*_X=T^*_{X/k}$ is trivial, i.e. $T^*_X\iso X\times V$, where $\dim V=r$. Then it holds that
\[\text{Chow}^n(T^*_{X_T/T}/{X_T})\iso \text{Chow}^n(V)\times X_T\]
and an analogous isomorphism for $U_T$.

Given a spectral datum $a\in \shf{A}_U(T)$, it is equivalent to a morphism $a: U_T\rightarrow \text{Chow}^n(V)$, which for ease of notation will also be denoted by $a$. By the key lemma \cite[Lemma 2.1]{CN}, it gives rise to a morphism
\[U_T\rightarrow V\times S^2V\times \cdots \times S^nV,\]
which in turn amounts to
\[U_T\rightarrow T^*_{X_T/T}\times S^2T^*_{X_T/T}\times\cdots \times S^nT^*_{X_T/T}\]
by the triviality of the relative cotangent bundle. This is exactly an element in
\[\bigoplus^n_{i=1} H^0(U_T, S^i T^*_{U_T/T})\]
\noindent \emph{Claim}: for each $i>0$, we have
\begin{align*}
H^0(X_T, S^i T^*_{X_T/T})\iso H^0(T, \sshf{T})\otimes_k H^0(X, S^i T^*_{X/k}).
\end{align*}
Given this \emph{Claim} and from Remark \ref{204}, we obtain that the natural map
\[H^0(X_T, S^i T^*_{X_T/T})\rightarrow H^0(U_T, S^i T^*_{U_T/T})\]
is an isomorphism for each $i>0$. Thus we get
\[\bar{a}\in \bigoplus^n_{i=1} H^0(X_T, S^i T^*_{X_T/T}).\]
We can view $\bar{a}$ as a morphism
\[ \bar{a}: X_T\rightarrow V\times S^2V\times \cdots \times S^nV,\]
which is equivalent to a section $\bar{a}\in \shf{A}_X(T)$.

To prove the \emph{Claim} above, consider the fibred product
$$\xymatrix{
X_T \ar[r]^{g}\ar[d]_{\varphi_T} &  X\ar[d]^{\varphi} \\
T\ar[r]^{h}&  \Spec(k).
}$$
It holds that
\[ T^*_{X_T/T}\iso g^* T^*_{X/k}.\]
Since $X$ is smooth, $T^*_{X/k}$ is locally free, hence
\[S^i T^*_{X_T/T}\iso g^*S^i( T^*_{X/k})\]
for all $i>0$. By flat base change, it therefore follows that
\begin{eqnarray*}
&&H^0(X_T, S^iT^*_{X_T/T})\\
&\iso & H^0(X_T, g^*S^i( T^*_{X/k}))\\
&\iso & H^0(\Spec{k}, h_*(\varphi_T)_*g^*S^i( T^*_{X/k}))\\
&\iso &H^0(\Spec{k}, h_*h^*\varphi_*S^i(T^*_{X/k}))\\
&\iso &H^0(\Spec{k}, h_*\sshf{T}\otimes_k \varphi_*S^i(T^*_{X/k}))\\
&\iso & H^0(T, \sshf{T})\otimes_k H^0(X, S^i(T^*_{X/k})). \end{eqnarray*}
This proves the \emph{Claim}.

In general, we take an affine open covering  $\{X_i\}$ of $X$ such that $T^*_{X/k}|_{X_i}$ is trivial for every $i$. Put $U_i=U\cap X_i$. For all $i$, it holds that
$\text{codim}(X_i-U_i, X_i)\ge 2$. By the preceding argument,
$\shf{A}_{X_i}(T)\rightarrow \shf{A}_{U_i}(T)$ is surjective for all $i$. Since $j^*$ is actually injective, we can lift section locally and then glue the resulting sections together. This finishes the proof of this lemma.
\end{proof}

\begin{thm}\label{303}
Let $\chi: X'\rightarrow X$ be a birational proper morphism of smooth varieties. Then the natural morphism $\shf{A}_X\rightarrow \shf{A}_{X'}$ induced by $\chi$ is an isomorphism as $k$-schemes.
\end{thm}
\begin{proof}
By the valuative criterion for properness, one can find an open subset $j: U\rightarrow X$ with $\text{codim}(X-U, X)\ge 2$ and a morphism $\nu : U\rightarrow X'$ such that $\chi \circ \nu=j$. As a result, we have the commutative diagram of schemes
$$\xymatrix{
\shf{A}_{X'} \ar@{^(->}[dr] &   \\
\shf{A}_X\ar[r]^{\iso} \ar@{^(->}[u] &  \shf{A}_{U},
}$$
where the bottom map is an isomorphism by Lemma \ref{302}.
Therefore $\shf{A}_X\rightarrow \shf{A}_{X'}$ is an isomorphism.
\end{proof}

\begin{exmp}
Let $Y \hookrightarrow X$ be a smooth subvariety. Denote by $\pi: \widetilde{X}={\rm BL}_Y(X)\rightarrow X$ the blow up of $X$ along $Y$. Then one has the injection of vector bundles
\begin{align*}
    0 \rightarrow \pi^* T^*_X \rightarrow T^*_{\widetilde{X}}.
\end{align*}
Note that ${\rm Chow}^n(T^*_X/X) \times_X \widetilde{X} \cong {\rm Chow}^n(\pi^* T^*_X/ \widetilde{X})$, and the injection above induces a morphism ${\rm Chow}^n(\pi^* T^*_X/ \widetilde{X}) \rightarrow {\rm Chow}^n(T^*_{\widetilde{X}} / \widetilde{X})$.

Now given a spectral datum $a_X: X \rightarrow {\rm Chow}_n(T^*_X /X)$, we define a new spectral datum $a_{\widetilde{X}}$ as the compositions of the following morphisms
\begin{align*}
    a_{\widetilde{X}}: \widetilde{X} \cong X \times_X \widetilde{X}   \xrightarrow{a_X \times 1} {\rm Chow}^n(T^*_X/X) \times_X \widetilde{X} \cong {\rm Chow}^n(\pi^* T^*_X/ \widetilde{X}) \rightarrow {\rm Chow}^n(T^*_{\widetilde{X}}/ \widetilde{X}),
\end{align*}
which gives a map $\shf{A}_X(k) \rightarrow \shf{A}_{\widetilde{X}}(k)$ between the $k$-points. By Theorem \ref{303}, the map is bijective.
\end{exmp}

\begin{cor}
For any projective ruled surface $X$, we have $\shf{A}_X \cong \shf{B}_X$. In particular it is an affine space.
\end{cor}
\begin{proof}
By Theorem \ref{303}, we can assume $X$ is minimal, i.e. it does not contain a smooth $(-1)$-rational curve. If $X=\prj{2}$, then $\shf{A}_X=\mathscr{B}_X$ is a point; otherwise $X\iso \mathbb{P}_C(E)$ for some smooth curve $C$ and a vector bundle $E$ of rank two, so $\shf{A}_X\iso \shf{A}_C$, and hence $\shf{A}_X=\shf{B}_X$, see \cite[Prop.~8.1]{CN2} and its following remark.
\end{proof}

\bigskip
\noindent\small{\textsc{School of Mathematics, Sun Yat-sen University}\\
W. 135 Xingang Rd., Guangzhou, Guangdong 510275, P.R.~China}\\
\emph{E-mail address}:  \texttt{songlei3@mail.sysu.edu.cn}

\bigskip
\noindent\small{\textsc{Department of Mathematics, South China University of Technology}\\
381 Wushan Rd., Guangzhou, Guangdong 510641, P.R.~China}\\
\emph{E-mail address}:  \texttt{hsun71275@scut.edu.cn}

\end{document}